\documentclass[12pt]{article}
\usepackage[russian]{babel}
\usepackage{amssymb,amsmath,amsthm}
\usepackage[cp1251]{inputenc}
\begin{document}

УДК 518.9
\begin{center}
\textbf{\LARGE Метод разрешающих функций для решения задачи преследования при
интегральных ограничениях на управления}\\[2mm]
{\bf Б.Т.Саматов}\\[2mm]
\end{center}

В настоящей работе основное внимание уделяется исследованию задачи
 преследования в линейных дифференциальных играх с интегральными
ограничениями. При исследовании этой задачи придерживаемся на
методе разрешающих функций.  Предложена иная конструкция
построения разрешающей функции, обосновывающая  правило
параллельного сближения игроков, т.е. П-стратегию для
преследователя. Для рассматриваемого случая получены новые
достаточные условия Л.С.Понтрягина к разрешимости задачи
преследования. В качестве примеров приводятся два класса игр,
каждый из которых может представить определенный интерес. Первый -
это контрольный пример Л.С.Понтрягина,  второй - задача простого
преследования при "$l$-поимки" убегающего. В этих
 предлагаемых примерах разрешающие функции играют ключевую
роль в решении задачи.
\begin{center}
\textbf{ 1. Введение}
\end{center}

Пусть в конечномерном евклидовом пространстве $\mathbb R^n$
движение объекта $z$ описывается линейным дифференциальным
уравнением
$$\dot{z}=A z+B u-C v, \eqno(1)$$
где $z\in R^{n}, \ u\in R^{m}, \ v\in R^k; \ n\geq 1, \ m\geq1, \
k\geq1; \ A, \ B, \ C$ - постоянные прямоугольные матрицы порядка
\ \ $n\times n, \ \ n\times m$ и \ $n\times k$ соответственно; $u$
- управляющий параметр преследования;\ $v$ - управляющий параметр
убегания. Параметры $u$ и $v$ выбираются в виде измеримых функций
$u=u(\cdot)$ и $v=v(\cdot)$ из класса $L_p[0,\ \infty)$,  $p>1$ и
удовлетворяют ограничениям
$$\int\limits^\infty_0|u(\tau)|^pd\tau \leq \rho^{p},\ \ \rho>0,  \eqno(2)$$
$$\int\limits^\infty_0|v(\tau)|^p d\tau \leq \sigma^{p}, \ \ \sigma>0. \eqno(3)$$
 Такие управления в дальнейшем будем называть допустимыми.

Терминальное множество  $M=M^0+M^1$, где $M^0$ - линейное
подпространство из $R^{n}$, \ а \ $M^1$ - выпуклое замкнутое
множество из ортогонального дополнения $L$ к подпространству $M^0$
в $R^{n}$.

{\sl Определение 1.}  В игре (1)-(3) из начального положения
$z_{0}\in R^{n}$ возможно завершение преследования за время
$T=T(z_{0})$, если по любой допустимой функции $v=v(t), \ 0\leq
t\leq T$, можно построить допустимую функцию $u(t, z_{0},
v_t(\cdot)) $, где $v_t(\cdot)=\{v(\tau), \ 0\leq\tau\leq t\},
0\leq t\leq T$, что  абсолютное непрерывное решение $z(t)$ задачи
Коши $\dot{z}=A z+B u(t,z_{0}, v_t(\cdot))-C v(t)$, $z(0)=z_{0}$,
попадает на терминальное множество $M$ за время, не превосходящее
числа $T=T(z_{0})$, т.е. $z(t^*)\in M$ при некотором $t^*\in[0,\
T]$. Число $T(z_{0})$ называется {\it гарантированным временем
преследования}.

Нахождение начальных положений (начальные состояния игроков,
параметры процесса и т.д.), из которых возможно завершение
преследования за конечное время, составляет задачу преследования.
Задаче преследования для различных классов игр посвящено много
работ и в зависимости от выбора стратегий в том или ином классе
были предложены ряд фундаментальных методов решения \ \ [1-5, 10,
11, 15, 18, 19, 23]. Основополагающими среди них являются работы
Л.С.Понтрягина [1-2] и  Н.Н.Красовского [3-4].

 В последние годы стали интенсивно исследоваться
дифференциальные игры  при наличии интегральных, разнотипных,
двойных(смешанных) и др. ограничений на управления игроков
[9,12,14,16,18,20,22,29-33,36-40]. В работах [12,14] были попытки
перенести, предложенный в [11] метод разрешающих функций для
дифференциальных игр преследования с геометрическими
ограничениями, на случай с интегральными ограничениями на
управления игроков.

Настоящая работа посвящается к  исследованию задачи преследования
для линейных дифференциальных игр с интегральными ограничениями.
При этом мы  придерживаемся на методе разрешающих функций[11] и
основываемся на формализации , предложенной Л.С.Понтрягиным в
[1-2], а так же используем некоторые идеи работ [9,16,19-22].
 В работе предложена иная конструкция построения разрешающей
 функции, чем в [12,14], обосновывающая правило параллельного сближения   при простых
движениях и интегральных ограничениях на управления игроков
[36-40]. Работа состоит из пяти пунктов. В первом пункте
приводится постановка задачи преследования при интегральных
ограничениях на управления игроков. Во втором,  для этой задачи
вводится и исследуется разрешающая функция $\lambda(t,\tau,v,z_0)$
и в третьем пункте  с помощью этой функции
 доказывается теорема о возможности завершения
преследования. В п.4 предлагаемый способ решения задачи
применяется к контрольному примеру Л.С.Понтрягина [1], [2].
Наконец, в п.5 достаточно подробно изучается задача преследования
при простом движении игроков для случая $l$-поимки. При этом для
преследователя предлагается стратегия, с помощью которой сближения
осуществляется наилучшим способом.

\begin{center}
\textbf{ 2. Построение разрешающей функции}

\end{center}

Пуст $\pi$-оператор ортогонального проектирования из $R^{n}$ на
подпространство $L$. Рассмотрим линейные отображения $\pi e^{A t}B
R^{m}\rightarrow L, \ \pi e^{A t}C R^k\rightarrow L$, когда $t\geq
0.$

{\sl Пpедположение 1.} Пусть существует непрерывная неособая
матрица $F(\cdot): R^k \rightarrow R^{m} $, являющаяся решением
матричного уравнения
$$\pi e^{A t}B X=\pi e^{A t}C,$$
где $X$- искомая матрица.

С помощью матрицы $F(\cdot)$ построим функцию(см.[15])
$$\chi^p(t)=\sup_{\int\limits_0^t|\omega(\tau)|^p d\tau\leq1} \int\limits_0^t|F(t-\tau)\omega(\tau)|^p d\tau,$$
где $ \omega(\cdot)-$ произвольная функция из пространства
$L_p[0,\infty)$ с указанным ограничением. С помощью функции
$\chi^{p}(t)$ определим величину
$$\nu^{p}=\sup_{0\leq t<\infty} \chi^{p}(t).$$
Очивидно, что $\nu^{p}>0.$

{\sl Пpедположение  2.} Пусть выполнено неравенство
$\rho^{p}>\sigma^{p}\nu^{p}.$

Если выполнено предположение 2, то очевидно, что  $\nu$ не может
равнятся $+\infty$.

 Введем многозначное отображение вида
$$U(t,\tau,v,\lambda)=(|F(t-\tau)v|^p+\lambda\delta)^\frac{1}{p}\pi e^{(t-\tau)A} B S-
\pi e^{(t-\tau)A}C v,\eqno(5)$$ где $0\leq\tau\leq t, \ v\in R^k,
 \delta=\rho^p-\sigma^p\nu^p,\ \lambda\geq0, \ S-$
шар радиуса единицы и с центром в нуле пространства $R^{m}$.

{\sl Лемма 1.} Включение $0\in U(t,\tau,v,\lambda)$ выполняется
при всех $ t, \tau, v, \lambda$, когда $0\leq\tau\leq t, \ v\in
R^k$ и $\lambda\geq0.$

{\sl Доказательство.} Рассмотрим многозначное отображение (5) при
$\lambda=0$. Тогда имеем $U(t,\tau,v,0)=|F(t-\tau)v|\pi
e^{(t-\tau)A} B S-\pi e^{(t-\tau)A}C v,$ где $0\leq\tau\leq t,\
v\in R^k. $  В силу предположения 1 получаем
$$\pi e^{(t-\tau)A} C v=\pi e^{(t-\tau)A} B F(t-\tau)v.$$
Отсюда находим, что
$$U(t,\tau,v,0)=|F(t-\tau)v| \pi e^{(t-\tau)A} B S-\pi e^{(t-\tau)A} BF(t-\tau)v=$$

$$=\left\{\begin{array}{cl} |F(t-\tau)v|\pi e^{(t-\tau)A}B\left(S-\frac{F( t-\tau) v}{|F(t-\tau)v|}\right), \ \mbox{ если
}F(t-\tau)v\neq 0,
\\[2mm]
\  \ \  \{0\}, \ \ \ \ \ \ \ \ \ \ \ \ \ \ \ \ \ \ \ \ \ \ \ \ \ \
\ \ \ \ \ \ \ \ \ \ \ \ \ \ \ \ \mbox{ если } F(t-\tau)v=0.
\\[2mm]\end{array}\right.
$$
Из того, что $F(t-\tau)v\in R^{m}$ и
$\frac{F(t-\tau)v}{|F(t-\tau)v|}\in S$, видно, что $0\in
U(t,\tau,v,0)$ при всех $0\leq\tau\leq t$ и $v\in R^k.$

Нетрудно проверить, что при $\delta>0$ многозначное отображение
(5) является монотонно возрастающим по включению по параметру
$\lambda\geq0,$ т.е. из $\lambda_1>\lambda_2$ следует включение
$U(t,\tau,v,\lambda_2)\subset U(t,\tau,v,\lambda_1).$
Cледовательно,  $0\in U(t,\tau,v,0)\subset U(t,\tau,v,\lambda)$
для всех $\lambda\geq0, \ 0\leq\tau\leq t$ и $v\in R^k$. Лемма
доказана.

{\sl Лемма 2.} Если $\pi e^{tA}z_0 \notin M^1$, где $z_0\in
R^{n}$, то функция
$$\lambda(t,\tau,v,z_0):=\max\left\{\lambda\geq0:\lambda(M^1-\pi e^{(t-\tau)A}z_0)
\cap U(t,\tau,v,\lambda)\neq\emptyset\right\}\eqno(6)$$ является
полунепрерывной сверху по переменным $\tau$ и $v$, где
$0\leq\tau\leq t$ и $v\in R^k.$

{\sl Доказательство.} Установим существование функции (6). Для
этого многозначные отображения $K(\lambda):=\lambda(M^1-\pi
e^{tA}z_0),$\ \ $U(\lambda):=U(t,\tau,v,\lambda)$ и
$Q(\lambda):=K(\lambda)\cap U(\lambda)$  будем рассматривать в
зависимости только от переменной $\lambda$, фиксируя произвольным
образом остальные переменные. Покажем, что область определения
многозначного отображения $Q(\lambda)$, т.е. $dom Q=\{\lambda:
Q(\lambda)\neq\emptyset\}$ (см.[8]), является непустым компактным
множеством.

Прежде покажем ограниченность $dom Q$. Для этого предположим
противное, т.е.  существует такое последовательность $\lambda_n\in
dom Q$, что при $n\rightarrow\infty$ получаем
$\lambda_n\rightarrow+\infty$. В силу теоремы I.I работы [8]
получаем, что $Q(\lambda)\neq\emptyset$ тогда и талько тогда,
когда выполняется неравенство
$$\min_{|\psi|=1}(W_{U(\lambda)}(\psi)+W_{K(\lambda)}(-\psi))\geq0,$$
где $\psi\in L$. Следовательно, из свойств опорных функций ([24])
и из конкретного вида многозначных отображений $K(\lambda)$ и
$U(\lambda)$ имеем
$$(|F(t-\tau)v|^p+\lambda\delta)^\frac{1}{p}\ W_{\pi e^{(t-\tau)A}B S} (\psi)-
(\pi e^{(t-\tau)A} C v, \psi)+\lambda W_{-\pi
e^{tA}z_0+M^1}(-\psi)\geq0\eqno(7)$$ для всех $\psi, \ |\psi|=1$.
Согласно предыдущей леммы
$$0\in (|F(t-\tau)v|^p+\lambda\delta)^\frac{1}{p} \pi e^{(t-\tau)A}B S-\pi e^{(t-\tau)A}C v.$$
Отсюда получаем, что
$$(|F(t-\tau)v|^p+\lambda\delta)^\frac{1}{p} W_{\pi e^{(t-\tau)A}B S}(\psi)-(\pi e^{(t-\tau)A}C v,\psi)\geq0$$
при всех $\psi, \ |\psi|=1$. Поэтому если $W_{-\pi
e^{tA}z_0+M^1}(-\psi)\geq0,$ то неравенство (7) выполняется для
всех $\lambda\geq0$. Остается рассмотреть случай, когда
 $W_{-\pi e^{tA}z_0+M^1}(-\psi)<0$. Из того, что $\pi e^{tA}z_0\notin M^1$ и
$M^1$-выпуклое замкнутое множество, получаем, что множество
$$\Gamma=\{\psi:|\psi|=1, \ W_{-\pi e^{tA}z_0+M^1}(-\psi)<0\}$$
является непустым множеством. Очевидно, что существует такое
$\psi$ из $\Gamma$ при котором неравенство (7) начиная с
некоторого $\lambda>0$ не выполняется, что и противоречить
предположению. Следовательно, множество $dom Q$ ограниченное.

Теперь остается показать замкнутость $dom Q$. Поскольку
многозначные отображения $U(\lambda)$ и $K(\lambda)$ для всех
$\lambda\geq0$ компактнозначны и непрерывны, то их опорные функции
$W_{U(\lambda)}(\psi)$ и $W_{K(\lambda)}(-\psi)$ так же будут
непрерывными для всех $\lambda\geq0$ и $\psi\in\Gamma$ (cм.[24]).
Тогда нетрудно проверить, что и функция
$\vartheta(\lambda):=\min_{|\psi|=1}[W_{U(\lambda)}(\psi)+W_{K(\lambda)}(-\psi)]$
является непрерывной по $\lambda, \ \lambda\geq0$. Отсюда и
следует замкнутость $dom Q=\{\lambda:\vartheta(\lambda)\geq0\}$,
что завершает доказательство компактности последнего.

Если $dom Q$-компактное множество из промежутка $[0,\infty)$, то
существует его наибольший элемент, которого и примем за функцию
$\lambda(t,\tau,v,z_0)$. Покажем полунепрерывность сверху этой
функции по переменным $\tau$ и $v$, где $0\leq\tau\leq t, \ v\in
R^k$.

Известно такой факт (см.[11]), что если некоторая функция $g(x,y)$
непрерывна на произведении компактов $X, \ Y$, где здесь $X, \
Y$-подмножества некоторых конечномерных евклидовых пространств, то
многозначное отображение
$$N(x)=\{y\in Y:g(x,y)\geq0\},$$
которое непусто для всех $x\in X$, является полунепрерывным сверху
на $X$.

Следовательно из непрерывности функции
$$\vartheta(\lambda,t,\tau,v,z_0)=\min_{|\psi|=1}\left[W_{U(t,\tau,v,\lambda)}(\psi)+
W_{\lambda(M^1-\pi e^{tA} z_0)}(-\psi)\right],$$ по  переменным
$\lambda$,  $\tau$ и $v$, где $\lambda\geq0$,   $0\leq\tau\leq t,
\ v\in R^k$, получаем, что многозначное отображение
$$Q(t,\tau,v,z_0)= dom Q =\{\lambda:\vartheta(\lambda,t,\tau,v,z_0)\geq0\}$$
будет полунепрерывным сверху по переменным $\tau$ и $v$. Отсюда и
из того, что $Q(t,\tau,v,z_0)(\subset[0,\infty))$-компактнозначное
отображение, легко получить, что и функция
$\lambda(t,\tau,v,z_0)=\max Q(t,\tau,v,z_0)$ является
полунепрерывным сверху по $\tau$ и $v$, где $0\leq\tau\leq t$ и
$v\in R^k$. Лемма доказана.

{\sl Определение  2.} Пусть выполнены предположения 1-2 и $\pi
e^{tA}z_0 \notin M^1$. Тогда в игре (1)-(3) {\it разрешающей
функцией} назовем функцию вида (6).

\begin{center}
\textbf{3. Теорема о возможности завершения преследования}
\end{center}

С помощью разрешающей  функции  $\lambda(t,\tau,v,z_0)$ вводится
функция
$$\Lambda(t,z_{0}):=1-\inf_{\int\limits_0^t|v(\tau)|^p d\tau
\leq\sigma^{p}}\int\limits^t_0\lambda(t,\tau,v(\tau),z_0)d\tau$$
Обозначим через $T=T(z_{0})$ первый положительный корень уравнения
$\Lambda(t,z_{0})=0$, если такого не существует то полагаем
$T(z_{0})=+\infty$.

{\sl Пpедположение 3.} Пусть для позиции  $z_{0}\in R^{n}$
существует конечный момент времени $T=T(z_{0})$.

{\sl Теорема .} Если выполнены предположения 1-2 и 3 , то из
позиции $z_{0}\in R^{n}$ в игре (1)-(3) возможно завершение
преследования за время $T=T(z_{0})$.

{\sl Доказательство.} Пусть для некоторого положения $z_{0}\in
R^{n}$ выполнены предположения 1-2 и 3. Рассмотрим функцию
$$\Lambda(t,z_{0},v(\cdot)):=1-\int\limits^t_0\lambda(T,\tau,v(\tau),z_{0})d\tau,$$
где $v=v(\cdot)$-произвольное управление убегающего, для которого
выполняется ограничение (3). Очевидно, что
$\Lambda(0,z_{0},v(\cdot))=1$, и  $\Lambda(t,z^0,v(\cdot))$
является равномерно непрерывной, монотонно невозрастающей функцией
по $t, \ 0\leq t\leq T$. Откуда в силу предположения 3
 следует существование такого момента $t^*$, что $t^*\leq
T$ и
$$\Lambda(t^*,z_{0},v(\cdot))=0.$$
При этом $\Lambda(t,z_{0},v(\cdot))>0$ для всех $t$, когда $0\leq
t<t^*$.

Теперь  рассмотрим многозначное отображение вида
$$M^1(\tau,v)=\{m^1\in M^1:\lambda(T,\tau,v,z_{0})(m^1-\pi e^{T A}z_{0})\in$$
$$\in (|F(t-\tau)v|^p+\lambda(T,\tau,v,z_{0})\delta)^\frac{1}{p}\pi e^{(T-\tau)A}B S-\pi e^{(T-\tau)A}C v\}. $$
Из того, что функция $\lambda(T,\tau,v,z_0)$ полунепрерывна сверху
по переменным $\tau$ и $v$ следует, что и многозначные отображения
$\lambda(T,\tau,v,z_0)(M^1-\pi e^{TA})z_{0}$ и
$(|F(t-\tau)v|^p+\lambda(T,\tau,v,z_{0})\delta)^\frac{1}{p}\pi
e^{(T-\tau)A}B S-\pi e^{(T-\tau)A}C v $  будут полунепрерывными
сверху по $\tau$ и $v$. Тогда в силу лемма 1.7.5 работы [25]
получаем, что многозначное отображение $M^1(\tau,v)$ является
измеримы. Следовательно, существует однозначная измеримая по
Борелю ветвь $m^1(\tau,v)\in M^1(\tau,v)$ (лемма 1.7.7. [25]).
Очевидно, что если функция $v=v(\tau)$ является суммируемой на
отрезке $[0, T]$, то и функция $m^1(\tau,v(\tau))$ является
суммируемой на $[0,T]$.  В ходе игры предпишем  преследователю
строить свой управление $u=u(\tau), \ 0\leq\tau\leq T,$ \
следующим образом:

 -до тех пор пока $\Lambda(t,z_{0},v(\cdot))>0$
управление $u=u(\tau)$ находится из уравнений
$$\pi e^{(T-\tau)A}B u(\tau)-\pi e^{(T-\tau)A}C v(\tau)=\lambda(T,\tau,v,z_{0})(m^1(\tau,v(\tau))-\pi e^{T A}z_{0}),$$
$$|u(\tau)|=(|F(T-\tau)v(\tau)|^p+\lambda(T,\tau,v(\tau),z_{0})\delta)^{\frac{1}{p}};\eqno(8)$$

-начиная с момента $t^*$ до $T(z_{0})$ управление $u=u(\tau), \
t^*\leq\tau\leq T$, находится из уравнений
$$\pi e^{(T-\tau)A}B u(\tau)=\pi e^{(T-\tau)A}C v(\tau),\ \ |u(\tau)|=(|F(T-\tau)v(\tau)|.\eqno(9)$$
Поскольку функции $\lambda(T,\tau,v(\tau),z_{0})$ и
$m^{1}(\tau,v(\tau))$ являются измеримыми по $\tau$, то по лемме
Филиппова А.Ф. [26] уравнения (8) и (9) разрешимы в классе
измеримых функций. Если это решение является неединственным, то
выбираем его наименьшим в лексикографическом смысле.

Покажем, что предложенный способ управления $u(\tau,v(\tau)), \
0\leq\tau\leq T$, позволяет завершить преследования при
произвольном допустимом управлении $v=v(\tau), \ 0\leq\tau\leq T$,
за время $T=T(z_{0})$. Для этого рассматривается следующая задача
Коши
$$\dot{z}=Az+Bu(\tau,v(\tau))-Cv(\tau),\ \ \  z(0)=z_{0}.$$
 Тогда  имеем
$$z(T)=e^{TA}z_{0}+\int\limits^{T}_0 e^{(T-\tau)A}[Bu(\tau,v(\tau))-Cv(\tau)]d\tau.$$
Из (8) и (9) находим
$$\pi z(T)=\pi e^{TA}z_{0}+\int\limits^{T}_0\pi e^{(T-\tau)A}[Bu(\tau,v(\tau))-C v(\tau)]d\tau=$$
$$=\pi e^{T A}z_{0}+\int\limits^{t^*}_0\lambda(T,\tau,v(\tau),z_{0})(m^1(\tau,v(\tau))-\pi e^{T A}z_{0})d\tau=$$
$$=\pi e^{TA}z_{0}(1-\int\limits^{t^*}_0\lambda(T,\tau,v(\tau),z_{0})d\tau)+
\int\limits^{t^*}_0\lambda(T,\tau,v(\tau),z_{0})m^1(\tau,v(\tau))d\tau.$$
Поскольку,  в момент $t^*$ выполняется равенство
$\Lambda(t^*,z_{0},v(\cdot))=0,$ а так же  в силу леммы  из [17]
находим, что

$$\pi z(T)=\int\limits^{t^*}_0\lambda(T,\tau,v(\tau),z_{0})m^1(\tau,v(\tau))d\tau\in
\int\limits^{t^*}_0\lambda(T,\tau,v(\tau),z_{0})M^1d\tau=$$
$$=M^1\int\limits^{t^*}_0\lambda(T,\tau,v(\tau),z_{0})d\tau=M^1,$$
 или $z(T)\in M$. Таким образом, игра (1) из начального положения
$z_{0}$  завершается за время $T=T(z_{0})$.

Осталось показать допустимость управления $u=u(\tau,v(\tau)), \
0\leq\tau\leq T$. Так как, по построению  выполняется соотношение
$$\int\limits^{T}_0|u(\tau,v(\tau))|^p d\tau=
\int\limits^{T}_0|F(T-\tau)v(\tau)|^p d\tau+
\delta\int\limits^{t^*}_0\lambda(T,\tau,v(\tau),z_{0})d\tau,$$ то
из легко проверяемого неравенства
$\int\limits^{T}_0|F(T-\tau)v(\tau)|^p d\tau\leq\sigma^p\ \nu^p$ \
\ \ и из того, что
$\int\limits^{t^*}_0\lambda(T,\tau,v(\tau),z_{0})d\tau=1,$
получаем $\int\limits^{T}_0|u(\tau,v(\tau))|^p
d\tau\leq\sigma^p\nu^p+\delta=\rho^p,$ так как
$\delta=\rho^p-\sigma^p\nu^p$.
 Теорема  доказана полностью.

{\sl Замечание}. Если для некоторых $z_{0}$ существует такой
момент $T=T(z_{0})$, что $\pi e^{TA}z_{0}\in M^1$ , то управление
$u=u(\tau,v(\tau)), \ 0\leq\tau\leq T$ находится только из
уравнений (9). Нетрудно показать, что в таком случае игра (1) из
точки $z_{0}$, так же завершается за время $T=T(z_{0}).$

\newpage

\begin{center}
\textbf{ 4. Контрольный пример Л.С.Понтрягина}
\end{center}

А. Пусть  движение  преследователя и убегающего задается
уравнениями
$$ \ddot{x}+\alpha \dot{x}=bu, \ x(0)=x_{0}, \  \dot{x}(0)=\dot{x}_{0},$$
$$ \ddot{y}+\beta\dot{y}=cv, \ \ y(0)=y_{0}, \  \dot{y}(0)=\dot{y}_{0},$$ где $x, y, v,  u\in R^n,  n\geq1, \
\alpha>0, \  \beta>0,  \ b>0, \ c>0 $;   управлении $u=u(\cdot) \
$ и $ \ v=v(\cdot)$- удовлетворяют ограничениям вида (2)-(3) при
$p=2$. Преследования считается завершенным, если в конечный момент
времени $x=y$.

Перейдем к соответствующей дифференциальной игре вида (1). Для
этого положим
$$(x-y, \ \dot{x}, \ \dot{y})=(z_{1}, \ z_{2}, \ z_{3})=z, \ \
(x_{0}-y_{0}, \ \dot{x}_{0}, \ \dot{y}_{0})=(z_{01}, \ z_{02}, \
z_{03})=z_{0},$$ и получаем эквивалентную систему уравнений
$$\dot{z}_{1}=z_{2}-z_{3},\
\dot{z}_{2}=-\alpha z_{2}+b u, \ \dot{z}_{3}=-\beta z_{3}+c v,\
 \eqno(10)$$
$$z_{1}(0)=z_{01},\ \ \ \   \ z_{2}(0)=z_{02}, \ \ \ \ \   z_{3}(0)=z_{03},$$
где $(z_{1}, z_{2}, z_{3})\in R^{3n}$  и $(z_{01}, z_{02},
z_{03})$ - начальное состояние игры (10)(см.[2]). Тогда
терминальное множество $M=\{(z_{1}, z_{2}, z_{3})\in
R^{3n}:z_{1}=0\}$. Отсюда $L=\{(z_{1},z_{2},z_{3})\in
R^{3n}:z_{2}=z_{3}=0\}.$ \ Матрицы $A, \ B$ и $C$ представимы в
виде

$$A=\left(\begin{array}{cl} \mathbb{O}\  \ \ \  \ \mathbb{E} \ \ \ \  \mathbb{E}
\\[2mm]
\mathbb{O} -\alpha \mathbb{E} \ \ \  \mathbb{O}
\\[2mm]
 \ \mathbb{O} \ \ \ \   \mathbb{O}  -\beta \mathbb{E}
\end{array}\right),
\ B=\left(\begin{array}{cl}   \mathbb{O}
\\[2mm]b
 \mathbb{E} \ \
\\[2mm]
\mathbb{O}
\end{array}\right), \
C=\left(\begin{array}{cl} \mathbb{O}
\\[2mm]
\mathbb{O}
\\[2mm]
-c \mathbb{E},
\end{array}\right),
$$
где $\mathbb{E}$ и $\mathbb{O}$-единичная и нулевая матрицы
порядка $n\times n$ соответственно.   Так как фундаментальная
матрица системы имеет вид

$$e^{tA}=\left(\begin{array}{cl} \mathbb{E} \ \
\frac{1-e^{-\alpha t}}{\alpha} \mathbb{E} \  -\frac{1-e^{-\beta t}}{\beta}\mathbb{E}\\[2mm]\mathbb{O} \ \ \ \ \  \ e^{-\alpha t} \mathbb{E} \ \ \ \ \ \ \ \ \ \  \mathbb{O}
\\[2mm]
 \ \mathbb{O} \ \ \ \ \ \ \ \ \   \mathbb{O} \  \ \ \ \ \  \ \  e^{-\beta t}\mathbb{E}
\end{array}\right),$$
то матричная функции $F(t), \ t\geq 0, $ удовлетворяющее
соотношению
$$\pi e^{A t}B F(t)=\pi e^{A t}C, $$
имеет вид $F(t)=\frac{c\alpha(1-e^{-\beta t})}{b\beta(1-e^{-\alpha
t})}\mathbb{E}, \ \ t\geq0$. Откуда находим, что
$$\chi^2(t)=\sup_{\int\limits_0^t|\omega(\tau)|^p d\tau\leq1}\int\limits^t_{0}f^2(t-\tau)|\omega(\tau)|^2d\tau,\ \  t\geq0,$$
где
$f(t-\tau)=\frac{c\alpha(1-e^{-\beta(t-\tau)})}{b\beta(1-e^{-\alpha(t-\tau)})},
\ 0\leq\tau\leq t.$ В силу непрерывности функции $f(t-\tau)$ по
$\tau, \ \ 0\leq\tau\leq t,$ нетрудно получить, что [15]
$$\nu^2=\sup_{0\leq t<\infty}\sup_{0\leq\tau\leq t}\frac{c^2\alpha^2(1-e^{-\beta(t-\tau)})}{b^2\beta^2(1
-e^{-\alpha(t-\tau)})}=\frac{c^2}{b^2}\max\{\frac{\alpha^2}{\beta^2},
\ \ 1\}.$$

\textbf{4.1. Разрешающая функция в Контрольном примере
Л.С.Понтрягина.} Полагаем, что выполнено предположение 2, т.е.
неравенство  $\rho>\sigma
  \max\{\frac{c\alpha}{b\beta}, \ \
1\}$. С учетом обозначений $\xi(t,z_0):=\pi
e^{At}z_0=z_{01}+\alpha (t)z_{02}-\beta (t)z_{03},$  где
$\alpha(t)=\frac{1-e^{-\alpha t}}{\alpha}, \ \
\beta(t)=\frac{1-e^{-\beta t}}{\beta}, \ t\geq0,$ вычислим
разрешающию функцию $\lambda=\lambda(t,\tau,v,z_0)$ для
рассматриваемого примера. Для этого предположим, что $z_0\neq0$, в
противно случае игры завершилась бы при $t=0$. Тогда из вида (5)
находим
$$U(t,\tau,v,\lambda)=[c^2 \ \beta^2(t-\tau)|v|^2+\lambda \ b^2 \ \alpha^2(t-\tau)\delta]^{\frac{1}{2}}S-c\beta(t-\tau)v,$$
где $0\leq\tau\leq t, \ \
\delta=\rho^2-\sigma^2\max\{\frac{c^2\alpha^2}{b^2\beta^2},1\}, \
\ v\in R^n$ и $S=\{u\in R^n:|u|\leq1\}$. Нетрудно показать, что
значения $\lambda$, для котрорых
$$-\lambda  \xi(t,z_0)+c \ \beta(t-\tau)v\in[c^2 \ \beta^2(t-\tau)|v|^2+
\lambda \ b^2 \ \alpha^2(t-\tau)\delta]^{\frac{1}{2}}S,$$
находятся в промежутке
$$0\leq\lambda\leq\frac{1}{|\xi(t,z_0)|^2}\max\{0, \ \delta \ b^2 \ \alpha^2(t-\tau)+2 \ c \ \beta(t-
\tau)(v,\xi(t,z_0))\}.$$ Откуда находим {\it разрешающию функцию
для Контрольного примера Л.С.Понтрягина}
$$\lambda(t,\tau,v,z_0)=\frac{1}{|\xi(t,z_0)|^2}\max\{0, \ \delta \ b^2 \ \alpha^2(t-\tau)+2 \ c \
 \beta(t-\tau)(v,\xi(t,z_0))\}.\eqno(11)$$
Из этой формулы видно, что функция $\lambda(t,\tau,v,z_0)$
является непрерывной по переменным $\tau$ и $v$, где
$0\leq\tau\leq t, \ \ v\in R^n$. Так как по предположению
$z_{01}\neq0$, то из конкретного вида $\xi(t,z_0)$ вытекает
существование такого промежутка $0\leq t\leq \bar{t}$, что
$\xi(t,z_0)\neq0$ на $[0,\bar{t}]$. Следовательно, функция
$\lambda(t,\tau,v,z_0)$ является непрерывной и по переменной $t$,
когда $0\leq t\leq \bar{t}$.

{\sl Лемма 4.} Если функция $\lambda(t.\tau,v,z_0)$ непрерывна по
совокупности переменных $t,\tau,v$, где $0\leq t\leq\bar{t}, \ \
0\leq\tau\leq t$ и $v\in R^n$, то и функция
$$\Lambda(t): = \inf_{\int\limits_0^t|v(\tau)|^2
d\tau\leq\sigma^{2}}\int\limits^t_0\lambda(t,\tau,v(\tau),z_{0})d\tau$$
- на промежутке $0\leq t\leq\bar{t}$ является непрерывной.

{\sl Доказательство.} Покажем непрерывность исследуемой функции в
точке $t_0=0$. Для этого зададим произвольное число
$\varepsilon>0$. Тогда в силу неравенства
$$|\inf_{\int\limits_0^t|v(\tau)|^2
d\tau\leq\sigma^{2}}\int\limits^t_0\lambda(t,\tau,v(\tau),z_{0})d\tau|\leq\sup_{\int\limits_0^t|v(\tau)|^2
d\tau\leq\sigma^{2}}|\int\limits^t_0\lambda(t,\tau,v(\tau),z_{0})d\tau|
$$ и непрерывности функции
$\int\limits^t_0\lambda(t,\tau,v(\tau),z_{0})d\tau$ по $t$, можно
указать такое $\delta=\delta(\varepsilon)>0$, что из $0\leq
t\leq\delta$ будет следовать неравенство
$|\Lambda(t)|<\varepsilon$.

Теперь покажем непрерывность функции $\Lambda(t)$ в произвольной
точке $t_0$, когда $0<t_0\leq\bar{t}$. Для этого надо показать,
что для любого $\varepsilon>0$ существует такой $\delta>0$, что из
$|t-t_0|<\delta$ следует неравенство
$|\Lambda(t)-\Lambda(t_o)|<\varepsilon$.

Нетрудно убедится, что
$$\Lambda(t_0)=\inf_{\int\limits_0^{t_{0}}|v(\tau)|^2
d\tau\leq\sigma^{2}}\int\limits^{t_0}_0\lambda(t_0,\tau,v(\tau),z_{0})d\tau=\inf_{\int\limits_0^t|v(\tau)|^2
d\tau\leq\sigma^{2}}\frac{t_0}{t}\int\limits^t_0\lambda(t_0,\frac{t_0}{t}\tau,
\ \sqrt{\frac{t}{t_0}}v(\tau),z_{0})d\tau.$$ В силу этого получаем
следующую последовательность соотношений
$$|\Lambda(t)-\Lambda(t_0)|=|\inf_{\int\limits_0^t|v(\tau)|^2
d\tau\leq\sigma^{2}}\int\limits^t_0\lambda(t,\tau,v(\tau),z_{0})d\tau-\inf_{\int\limits_0^t|v(\tau)|^2
d\tau\leq\sigma^{2}}\frac{t_0}{t}\int\limits^t_0\lambda(t_0,\frac{t_0}{t}\tau,
\ \sqrt{\frac{t}{t_0}}v(\tau),z_{0})d\tau|$$
$$\leq \sup_{\int\limits_0^t|v(\tau)|^2
d\tau\leq\sigma^{2}}\left|\int\limits^t_0\lambda(t,\tau,v(\tau),z_{0})d\tau-
\frac{t_0}{t}\int\limits^t_0\lambda(t_0,\frac{t_0}{t}\tau, \
\sqrt{\frac{t}{t_0}}v(\tau),z_{0})d\tau\right|\leq$$
$$\leq \sup_{\int\limits_0^t|v(\tau)|^2
d\tau\leq\sigma^{2}}\int\limits^t_0|\lambda(t,\tau,v(\tau),z_{0})d\tau-\frac{t_0}{t}\lambda(t_0,\frac{t_0}{t}\tau,
\ \sqrt{\frac{t}{t_0}}v(\tau),z_{0})|d\tau\leq$$
$$\leq \sup_{\int\limits_0^t|v(\tau)|^2
d\tau\leq\sigma^{2}}\int\limits^t_0[|\lambda(t,\tau,v(\tau),z_{0})-\lambda(t_0,\tau,v(\tau),z_{0})|+
|\lambda(t_0,\tau,v(\tau),z_{0})-\lambda(t_0,\frac{t_0}{t}\tau, \
\sqrt{\frac{t}{t_0}}v(\tau),z_{0})|+$$
$$+|\lambda(t_0,\frac{t_0}{t}\tau, \ \sqrt{\frac{t}{t_0}}v(\tau))-
\frac{t_0}{t}\lambda(t_0,\frac{t_0}{t}\tau,\sqrt{\frac{t}{t_0}}v(\tau),z_{0})|]d\tau.$$
Когда $|t-t_0|<\delta$, то существует такое число
$\bar{\delta}=\bar{\delta}(\delta)>0$, что
$|\frac{t_0}{t}\tau-\tau|<\bar{\delta}$ и
$|\sqrt{\frac{t}{t_0}}v-v|<\bar{\delta}$. Отсюда, согласно
непрерывности функции $\lambda(t,\tau,v,z_{0})$ в каждой точке
$(t_0,\tau,v)$, получаем, что \ \
$|\lambda(t,\tau,v,z_{0})-\lambda(t_0,\tau,v,z_{0})|<\varepsilon_1,
\ \ \ \ \ |\lambda(t_0,\tau,v)-\ \ \ \ \ \ \ \ \
\lambda(t_0,\frac{t_0}{t}\tau, \
\sqrt{\frac{t}{t_0}}v)|<\varepsilon_2,\ \
|\lambda(t_0,\frac{t_0}{t}\tau, \
\sqrt{\frac{t}{t_0}}v)|\cdot\frac{|t-t_0|}{t}<\varepsilon_3.$
Отсюда находим
$$|\Lambda(t)-\Lambda(t_0)|<t(\varepsilon_1+\varepsilon_2+\varepsilon_3)\leq\bar{t}(\varepsilon_1+\varepsilon_2+\varepsilon_3)<\varepsilon,$$
что и требовалось доказать.

\textbf{4.2. Решение задачи преследования.}  В дальнейшем
рассматривается уравнение
$$1-\inf_{\int\limits_0^{t}|v(\tau)|^2
d\tau\leq\sigma^{2}}\int\limits^t_{0}\lambda(t,\tau,v(\tau),z_0)d\tau=0.\eqno(12)$$

{\sl Утверждение 1.}  Если
$\rho^2>\sigma^2\frac{c^2}{b^2}\max\{\frac{\alpha^2}{\beta^2}, \ \
1\}$, то для каждого $z_0=(z_{01},z_{02},z_{03})\in R^{3n}$ при
$z_{01}\neq0$, существует конечное решение $t=T(z_0)$ уравнения
(12).

{\sl Доказательство.}  Из очевидного неравенства
$\int\limits^t_{0}\max \{0,\varphi(\tau)\}d\tau\geq\max
\{0,\int\limits^t_{0}\varphi(\tau)d\tau\},$ где $\varphi(\tau)$-
произвольная суммируемая функция, получаем, что
$$J(t,v(\cdot),z_0)=\int\limits^t_{0}\max\{0,\delta \ b^2 \ \alpha^2(t-\tau)+2c \ \beta(t-\tau)(v(\tau),
\ \xi(t,z_0))\}d\tau\geq$$
$$\geq\{0,\delta \ b^2\int\limits^t_{0}\alpha^2(t-\tau)d\tau+2\int\limits^t_{0}c \ \beta(t-\tau)(v(\tau),
\ \xi(t,z_0))\}d\tau.$$ С силу неравенства Коши-Буняковского и из
определения функции $\chi^2(t)$ имеем
$$\int\limits^t_{0}\left(\frac{c \ \beta(t-\tau)}{b \
\alpha(t-\tau)}v(\tau), \ b \ \alpha(t-\tau)\xi(t,z_0)\right)
d\tau\leq b
\left(\int\limits^t_{0}\frac{c^2\beta^2(t-\tau)}{b^2\alpha^2(t-
\tau)}|v(\tau)|^2d\tau\right)^{\frac{1}{2}}\times$$
$$\times\left(\int\limits^t_{0}\alpha^2(t-\tau)d\tau\right)^\frac{1}{2}
|\xi(t,z_0)|\leq b\sigma
\nu\left(\int\limits^t_{0}\alpha^2(t-\tau)d\tau\right)^\frac{1}{2}|\xi(t,z_0)|.$$
Следовательно,
$$J(t,v(\cdot),
z_0)\geq\max\{0, \ \delta \
b^2\int\limits^t_{0}\alpha^2(t-\tau)d\tau-2\sigma\nu
b\left(\int\limits^t_{0}\alpha^2(t-\tau)d\tau\right)^\frac{1}{2}|\xi(t,z_0)|\}.$$
Из того, что $\delta=\rho^2-\nu^2\sigma^2>0$, функция $\xi(t,z_0)$
при $\alpha>0$ и $\beta>0$ ограниченная на $0\leq t<+\infty$, и
$\int\limits^t_{0}\alpha^2(t-\tau)d\tau$ монотонно возрастающая по
$t$ при $t\geq 0$; получаем существование такого момента
$\theta=\theta(z_0)$ для которого выполняется равенство
$$\delta\ b^2\int\limits^{\theta}_{0}\alpha^2(\theta-\tau)d\tau-2\sigma\nu b\left(\int\limits^{\theta}_{0}\alpha^2(\theta-
\tau)d\tau\right)^\frac{1}{2}|\xi(\theta,z_0)|=|\xi(\theta,z_0)|^2$$
Из этого равенства получаем, что
$$|\xi(\theta,z_0)|=(\rho-\sigma\nu)b\left(\int\limits^{\theta}_{0}\alpha^2(\theta-\tau)d\tau\right)^\frac{1}{2}.$$
Следовательно, для произвольного допустимого управления
$v=v(\tau)$ при $0\leq\tau\leq \theta$ выполняется соотношение
$$\inf_{v(\cdot)}J(v(\cdot), \theta, \ z_0)-|\xi(\theta,z_0)|^2\geq0.$$
Откуда имеем
$$\Lambda^\ast(\theta):=1-\Lambda(\theta)=1-\inf_{\int\limits_0^\theta|v(\tau)|^2
d\tau\leq\sigma^{2}}\int\limits^{\theta}_{0}\lambda(\theta,\tau,v(\tau),z_0)d\tau\leq0.$$
Из этого и в силу леммы 4  получаем, что уравнение (12) имеет для
всех $z_0\notin M$ конечное решение $T(z_0)\leq \theta.$
Утверждение 1 доказано.

{\sl Утверждение 2.}   Если
$\rho^2>\sigma^2\frac{c^2}{b^2}\max\{\frac{\alpha^2}{\beta^2}, \ \
1\}$, то в игре (10) для каждого $z_0=(z_{01},z_{02},z_{03})\in
R^{3n}$, когда $z_{01}\neq0$, при помощи разрешающей функции (11)
возможно завершение преследования за время $T(z_{0})$, где
$T(z_{0})$-первый положительный корень уравнения (12).

{\sl Доказательство.} Пусть убегающий выбирает
 произвольное допустимое управление $v=v(\tau)$, $0\leq \tau \leq T(z_{0})$, а  преследователя с
 строить своё управление следующим образом:

$$u(\tau,v(\tau))=\left\{\begin{array}{cl}
\frac{c\beta(T-\tau)}{b\alpha(T-\tau)}v(\tau)-\lambda^{*}(T,\tau,v(\tau),z_0)\xi(T,z_0),\mbox{
когда } 0\leq\tau\leq t^*
\\[2mm]
\frac{c\beta(T-\tau)}{b\alpha(T-\tau)}v(\tau), \ \ \ \ \ \ \ \ \ \
\ \ \ \ \ \ \ \ \ \ \  \ \ \ \ \ \ \ \ \ \ \ \ \mbox { когда }
t^*<\tau\leq T,
\\[2mm]\end{array}\right.\eqno(13)
$$
где
$$\lambda^{*}(T,\tau,v(\tau),z_0)=\frac{1}{|\xi(T,z_0)|^2}\max\{0, \ \delta \ b \ \alpha(T-\tau)
+2\frac{c\beta(T-\tau)}{b\alpha(T-\tau)}(v(\tau),\xi(T,z_0))\},$$
а $t^*$-первый положительный корень уравнения
$$\int\limits^t_{0}\max\{0, \ \delta \ b^2 \
\alpha^2(T-\tau)+2c\beta(T-\tau)(v(\tau),
\xi(T,z_0))\}d\tau=|\xi(T,z_0)|^2.\eqno(14)$$ Существование такого
корня следует из утверждения 1.

Прежде покажем допустимость предложенного управления (13). Для
этого из вида функции (13) нетрудно вычислить, что
$$|u(\tau,v(\tau))|^2=\left\{\begin{array}{cl}
\frac{c^2\beta^2(T-\tau)}{b^2\alpha^2(T-\tau)}|v(\tau)|^2+\frac{\delta}{|\xi(T,z_0)|^2}\max\{0,
\ \delta \
b^2\alpha^2(T-\tau)+2c\beta(T-\tau)(v(\tau),\xi(T,z_0))\},
\\[2mm]
\ \ \ \ \ \ \ \ \ \ \ \ \ \ \ \ \ \ \ \ \ \ \ \ \ \ \ \ \ \ \ \ \
\ \ \mbox{ при } 0\leq\tau\leq t^*,
\\[2mm]
\frac{c^2\beta^2(T-\tau)}{b^2\alpha^2(T-\tau)}|v(\tau)|^2, \ \ \ \
\ \ \ \ \ \ \ \ \ \ \mbox{ при }t^*<\tau\leq T,
\\[2mm]\end{array}\right.
$$
и проинтегрируем ее в промежутке от $0$ до $T$:
$$\int\limits^{T}_{0}|u(\tau,v(\tau))|^2d\tau=\int\limits^{T}_{0}\frac{c^2\beta^2(T-\tau)}{b^2\alpha^2(T-
\tau)}|v(\tau)|^2d\tau+\frac{\delta}{|\xi(T,z_0)|^2}\times$$
$$\times\int\limits^{t^*}_{0}\max\{0, \ \delta \ b^2\alpha^2(T-\tau)+2c\beta(T-\tau)(v(\tau),\xi(T,z_0))\}d\tau.$$
В силу того, что $t^*$-решение уравнения (14), и из определения
функция $\chi^2(t), \ 0\leq t<\infty,$ а так же из неравенства
$\chi^2(t)\leq\nu^2$ находим, что
$$\int\limits^{T}_{0}|u(\tau,v(\tau))|^2d\tau=\int\limits^{T}_{0}\frac{c^2\beta^2(T-\tau)}{b^2\alpha^2(T-
\tau)}|v(\tau)|^2d\tau+\delta\leq\sigma^2\chi^2(t)+\delta\leq\sigma^2\nu^2+\delta.$$
Так как $\delta=\rho^2-\sigma^2\nu^2$, то выполнено ограничение
$\int\limits^{T}_{0}|u(\tau,v(\tau))|^2d\tau\leq\rho^2,$ что и
требовалось показать.

Теперь покажем завершение игры в промежутке времени
$[0,T(z_{0})]$, при применении со стороны преследователя
управлению (13), а убегающий произвольное допустимое управление
$v=v(\tau), \ 0\leq\tau\leq T$.  Подставив эти управления в
уравнения (10)  по формуле Коши получаем
$$\pi z(T)=\pi e^{AT}z_0+\int\limits^{T}_{0}\pi e^{A(T-\tau)}(B u(\tau,v(\tau))-Cv(\tau))d\tau=$$
$$=\xi(T,z_0)+\int\limits^{T}_{0}(b\alpha(T-\tau)u(\tau,v(\tau))-c\beta(T-\tau)v(\tau))d\tau=\xi(T,z_{0})-$$
$$-\frac{\xi(T,z_0)}{|\xi(T,z_0)|^2}\int\limits^{t^*}_{0}\max\{0, \ \delta \ b^2\alpha^2(T-\tau)+2c\beta(T
-\tau)(v(\tau), \
\xi(T,z_0))\}d\tau=\frac{\xi(T,z_0)}{|\xi(T,z_0)|}$$
$$\left(|\xi(T,z_0)|^2-\int\limits^{t^*}_{0}\max{0,\delta b^2\alpha^2(T-\tau)+2c\beta(T-\tau)
(v(\tau),\xi(T,z_0))}d\tau\right)=0$$ Это означает, что $\pi
z(T)=z_{01}(T)\in M$ или $x(T)=y(T).$ Таким образом завершается и
доказательство утверждении 2.

\begin{center}
\textbf{5. Преследование при простом движении игроков. Случай
$l$-поимки}
\end{center}

В предыдущем параграфе мы рассмотрели пример, который имеет
достаточно общий характер. В настоящем пункте рассматривается
дифференциальная игра при простом движении игроков для случая
$l$-поимки. Предлагается разрешающая функция, при помощи которой
сближение игроков происходит наилучшим образом. Это игра для
случая геометрических ограничений исследованы в работах [7,
11,15,21,28] и в др.

Пусть в пространстве $R^n$ точка $x$ со скоростью $u$ преследует
точку $y$, которая  перемещается в пространстве $R^n$ со скоростью
$v$. Их движения описываются уравнениями
$$\dot{x}=u, \ \ x(0)=x_{0},\ \  \ \dot{y}=v, \ \ y(0)=y_{0}, \eqno(15)$$
где $x_{0}, y_{0}$- начальные положения точек;  скорости  $u$ и
$v$ выбираются в виде измеримых функций
$u(\cdot):[0,\infty)\rightarrow R^n$ и $\
v(\cdot):[0,\infty)\rightarrow R^n$ из пространства
$L_2[0,\infty)$ и удовлетворяют ограничениям
$$\int\limits^\infty_0|u(\tau)|^2d\tau\leq\rho^2, \ \rho>0, \eqno(16)$$
$$\int\limits^\infty_0|v(\tau)|^2d\tau\leq\sigma^2, \ \sigma\geq0.\eqno(17)$$

 В дальнейшем,  точку $x$ назовем преследователем, а точку $y$ убегающим  и  процесс
 преследование считаем завершенным, если в некоторый
конечный момент времени  выполнено соотношение $|x-y|\leq l,$ где
$l\geq0$. В начальный момент времени $t=0$ полагаем, что
$|x_{0}-y_{0}|>l$.

Для удобства исследования вводится переменное $z=x-y$. Тогда
уравнения (15) относительно  $z$ принимает вид
$$\dot{z}=u-v, \ z(0)=z_{0}, \eqno(18)$$
где $z_{0}=x_{0}-y_{0}, \ \ \ |z_{0}|>l$;\  терминальное множество
$M$ представляется в виде $M=\{z\in R^n:|z|\leq l\}=lS,\ \
l\geq0,$ здесь $S$-шар радиуса единицы с центром в нуле
пространства $R^n$.

Для конструирование стратегии преследователю  разрешается
использовать в каждый  момент времени только текущее значение
управления $v(t)$ и постоянные $z_0, \ \rho, \ \sigma$ и $l$.

\textbf{5.1. Разрешающая функция в случае $l$-поимки}.  С помощью
более общего метода  предложенного в п.2 построим разрешающую
функцию и здесь. Поскольку, для игры (18) предположение 2
принимает вид $\rho>\sigma,$ то из  (6) имеем функцию
$$\lambda(v,z_{0})=\max\{\lambda\geq0: \lambda(M-z_{0})\cap
U(\lambda,v)\neq\emptyset \} ,$$ где $U(\lambda,
v)=(|v|^2+\lambda\delta)^{1/2} S-v,  \ \delta=\rho^2-\sigma^2.$
Найдем те $\lambda\geq0$, которые удовлетворяют соотношению
$$\lambda(lS-z_0)\cap U(\lambda,v)\neq\emptyset.$$
Это соотношение эквивалентно неравенству
$$\lambda W_{lS-z_0}(-\psi)+W_{U(\lambda,v)}(\psi)\geq0,$$
для всех $\psi\in R^n$, при $|\psi|=1$. Отсюда
$$\lambda((\psi,z_0)+l|\psi|)+(|v|^2+\lambda\delta)^{1/2}|\psi|-(\psi,v)\geq0,$$
или
$$\lambda l+(|v|^2+\lambda\delta)^{1/2}\geq(\psi,v-\lambda z_0).$$
Однако $\max_{|\psi|=1}(\psi,v-\lambda z_0)=|v-\lambda z_0|.$
Тогда получаем $\lambda
l+(|v|^2+\lambda\delta)^{1/2}\geq|v-\lambda z_0|.$ Выполнив
элементарные вкладки, находим, что
$$0\geq\lambda
[\lambda^2h^2-2\lambda(h(\delta+2(v,z_0))+2l^2\delta)+(\delta+2(v,z_0))^2-4l^2|v|^2
],$$ где $h=|z_0|^2-l^2>0.$  Отсюда получаем, что
$0\leq\lambda\leq\max\{0,\lambda^*(v,z_0)\}),$ где
$$\lambda^*(v,z_0)=\frac{1}{h^2}(h(\delta+2(v,z_0))+2l^2\delta+|z_0\delta+vh|2l).$$
Таким образом, в этой задаче  разрешающую функцию $\lambda(v,z_0)$
определяем в виде
$$\lambda(v,z_0):=max\{0,\lambda^*(v,z_0)\}.\eqno(19)$$ Нетрудно проверить,
что

$$\lambda(v,z_0)=\left\{\begin{array}{cl} >0, \ \ \ \mbox {если} \ \ \ \delta+2(v,z_0)+2l|v|>0,\\[2mm]
0, \ \ \ \ \  \ \mbox {если} \ \ \ \delta+2(v,z_0)+2l|v|\leq0.
\end{array}\right.$$

\textbf{5.2. Построение стратегии}. Пусть в каждый момент времени
$t\geq0$ преследователю  известна информация только о текущих
значениях управления $v(t)$ и начальных заданных $z_0, \ \rho, \
\sigma$ и $l$. Тогда с помощью разрешающей функции
$\lambda(v,z_0)$ построим стратегию. Для этого аналогично
предложенному в пункте 3 способу решения задачи преследования
рассмотрим многозначное отображение $$M(v,z_{0})=\{m\in lS:
\lambda(v,z_0)(m-z_{0})\in(|v|^2+\lambda(v,z_{0})\delta)^{1/2}S-v\}.$$
В силу определения функции $\lambda(v,z_0)$ легко заметить, что
существует единственная однозначная измеримая по Борелю ветвь
$m(v,z_{0})\in M(v,z_{0})$. Из вида $M(v,z_{0})$, используя
свойства опорных функций,  вычислим  функцию $m(v,z_{0})$ при
помощи равенства
 $$\min_{|\psi|=1}[\lambda(v,z_0)W_{lS-z_{0}}(\psi)+W_{(|v|^2+\lambda(v,z_{0})\delta)^{1/2}S-v}(-\psi)]=0.$$
 Откуда получаем
 $$\min_{|\psi|=1}[\lambda(v,z_0)(|\psi|l-(\psi,z_{0}))+|\psi|(|v|^2+\lambda(v,z_{0})\delta)^{1/2}+(\psi,v)]=0,$$
или
$$\lambda(v,z_0)l+(|v|^2+\lambda(v,z_{0})\delta)^{1/2}+\min_{|\psi|=1}(\psi,v-\lambda(v,z_{0})z_{0})=0.$$
Из явного вида $\lambda(v,z_{0})$ нетрудно проверить, что
$\lambda(v,z_0)l+(|v|^2+\lambda(v,z_{0})\delta)^{1/2}=|v-\lambda(v,z_{0})z_{0}|.$
Следовательно получаем соотношение
$$\min_{|\psi|=1}(\psi,v-\lambda(v,z_{0})z_{0})=-|v-\lambda(v,z_{0})z_{0}|,$$
которое выполняется при
$\psi^{*}:=-(v-\lambda(v,z_{0})z_{0})/|v-\lambda(v,z_{0})z_{0}|.$
Отсюда и получаем явный вид функции:  $m(v,z_{0})=l\psi^{*}$. В
силу равенства (8) определим стратегию для преследователя.

{\sl Определение  3.} Пусть $\rho>\sigma$.Тогда в игре (15)-(18)
$\Pi_{l}-$ стратегией преследователя назовем функцию вида

$$u(v,z_{0})=v+\lambda(v,z_{0})(m(v,z_{0})-z_{0}),\eqno(20)$$
где
$m(v,z_{0})=-\frac{v-\lambda(v,z_{0})z_{0}}{|v-\lambda(v,z_{0})z_{0}|}l,
\ \lambda(v,z_{0})$- разрешающая функция вида (19).

Из вида стратегии (20) легко вычислить, что
$$|u(v,z_{0})|^2=|v|^2+\delta\lambda(v,z_{0}).\eqno(21)$$

\textbf{5.2.О возможности $l$-поимки}. Далее, рассматривается
уравнение
$$1-\int\limits^t_0\lambda(v(\tau),z_{0})d\tau=0\eqno(22)$$
Обозначим через $T=T(z_{0},v(\cdot))$ первый положительный корень
этого уравнения, где $v(\cdot)$-произвольное допустимое управление
убегающего. Существование и ограниченность такого корня
показывается в следующем утверждении.

{\sl Утверждение 3.} . В игре (15)-(18) из точки $z_{0}\notin lS$
при помощи  $\Pi_{l}-$ стратегии (20) возможно завершение
преследования за время $T(z^0,v(\cdot))$, удовлетворяющее
неравенству
$$T(z_{0},v(\cdot))\leq
\left(\frac{|z_{0}|-l}{\rho-\sigma}\right)^2,$$ при всех
допустимых управлений  $v(\cdot).$

Доказательство. В начале покажем ограниченность и существование
корня $T(z_{0},v(\cdot))$ уравнения (21). Для этого учитывая
ограничение на управления $v=v(t), \ t\geq0$, в виде (17), и
используя неравенство Коши-Буняковского, получаем
$$1-\int\limits^t_0\lambda(v(\tau),z_{0})d\tau=1-\int\limits^t_0\frac{1}{h^2}\max\{0, h(\delta+2(v(\tau),z_{0}))+$$
$$+2l^2\delta+2l|z_{0}\delta+v(\tau)h|\}d\tau \leq$$
$$\leq1-\frac{1}{h^2}\max\{0, \ t\delta(h+2l^{2})+2h\int\limits^t_0(v(\tau),z_{0})d\tau+2l\int\limits^t_0|z_{0}\delta
+v(\tau)h|d\tau\})\leq$$
$$\leq1-\frac{1}{h^2}\max\{0, \ t\delta(h+2l^2)-2h\int\limits^t_0|v(\tau)||z_{0}|d\tau+2l\int\limits^t_0(|z_{0}|
\delta-|v(\tau)|h)d\tau=$$
$$=1-\frac{1}{h^2}\max\{0, t\delta(h+2l^2+2l|z_{0}|)-(2h|z_{0}|+2lh)\int\limits^t_0|v(\tau)|d\tau\}\leq$$
$$\leq 1-\frac{1}{h^2}\max\{0, t\delta(h+2l^2+2l|z_{0}|)-2h(|z_{0}|+l)\sqrt{t}(\int\limits^t_0
|v(\tau)|^2d\tau)^\frac{1}{2}\leq$$
$$\leq 1-\frac{1}{h^2}\max\{0, t\delta(|z_{0}|+l)^2-2h\sqrt{t}\sigma(|z_{0}|+l)\}=$$
$$=1-\frac{1}{h^2}\max\{0, t(\rho^2-\sigma^2)(|z_{0}|+l)^2-2h\sqrt{t}\sigma(|z_{0}|+l)\}.$$
Нетрудно проверить, что последняя часть в этой цепочке неравенств
обрашается в нуль в момент времени $t=\theta$, где
$$\theta=\left(\frac{|z_{0}|-l}{\rho-\sigma}\right)^2.$$
Следовательно, в силу непрерывности по $t, t\geq0$ функции
$$\Lambda(t,v(\cdot))=1-\int\limits^t_0\lambda(v(\tau),z_{0})d\tau,$$
 имеем $\Lambda(\theta,v(\cdot))\leq0$. Отсюда и из того что
$\Lambda(0,v(\cdot))=1$ находим существование такого момента
$T=T(z^0,v(\cdot))$, что $\Lambda(T,v(\cdot))=0$. При этом
очевидно, что $T(z^0,v(\cdot))\leq\theta$, что и требовалось
показать.

Теперь покажем, что из произвольной точки $z_{0}\notin lS$
преследование завершается именно в этот момент времени
$T(z^0,v(\cdot))$. Для этого при произвольном допустимом
управлении убегающего  $v=v(t), \ 0\leq t\leq T$, преследователю
предпишем реализацию стратегии (20) т.е. управлении
$u=u(v(t),z_{0}), \ 0\leq t\leq T $. Тогда уравнение (18)
приводится к виду
$$\dot{z}=\lambda(v(t),z_{0})(m(v(t),z_{0})-z_{0}),$$
$$z(0)=z_{0},  \ 0\leq t\leq T.$$
Отсюда по формуле Коши получаем
$$z(t)=z_{0}+\int\limits^t_0\lambda(v(\tau), \ z_{0})(m(v(\tau),z_{0})-z_{0})d\tau,$$
где $0\leq t\leq T$. Учитывая соотношение
$1-\int\limits^T_0\lambda(v(\tau),z_{0})d\tau=0$ и вид функции
$m(v(\tau),z_{0})$ имеем
$$|z(T)|-l=|z_{0}-z_{0}\int\limits^T_0\lambda(v(\tau),z_{0})d\tau+\int\limits^T_0\lambda(v(\tau),z_{0})
m(v(\tau),z_{0})d\tau|-l\leq$$
$$\leq|z_{0}|(1-\int\limits^T_0\lambda(v(\tau),z_{0})d\tau)-l(1-\int\limits^T_0\lambda(v(\tau),z_{})d\tau)
=0,$$ или $|z(T)|\leq l$, что и требовалось показать.

Допустимость реализации $\Pi_{l}-$стратегии (20) следует из (21) и
(17), и из того,что $T$-есть решение уравнения (22):
$$\int\limits^T_0|u(v(\tau),z_{0})|^2d\tau=\int\limits^T_0|v(\tau)|^2d\tau+
\delta\int\limits^T_0\lambda(v(\tau),z_{0})d\tau\leq\sigma^2+
\delta\int\limits^T_0\lambda(v(\tau),z_{0})d\tau=\sigma^2+\delta=\rho^2,$$
 что  завершает доказательство утверждение 3.

{\sl Примечание 1.} В игре (15)-(17), если  $\rho>\sigma$ и $l=0$,
то разрешающая функция принимает простой вид
$\lambda(z_{0},\delta,v):=
\frac{1}{|z_{0}|^{2}}max\{0,\delta+2(v,z_{0})\}.$ В работе [38]
стратегия параллельного сближения:
$u(v,\delta,z_{0})=v-\lambda(v,z_{0})z_{0}$  аналогично к работе
[15], определяется как $\Pi_{I}-$ стратегия ( $I$- означает , что
ограничения на управления игроков интегральные). Показано, что при
реализации $\Pi_{I}-$ стратегии, преследователь не выпускает
убегающего из замкнутого шара $S(x_{0},y_{0},\rho,\sigma): =
\{w\in R^{n}: \ |w-y_{0}+\sigma^{2}z_{0}/\delta|\leq |z_{0}|\rho
\sigma/\delta\},$
 т.е. в силу утверждения 3, если в некоторый момент времени $T$ осуществляется равенство
$x(T)=y(T)$, то $y(T)\in S(x_{0},y_{0},\rho,\sigma)$, при
произвольном допустимом управлении убегающего. В случае
пересечении этого шара с некоторым заранее заданным множеством  (
с "линией жизни" [5]),  можно показать возможность попадании
убегающего в некоторую точку из этого пересечения и  до момента
попадания оставаться непойманным.\\

Пользуясь случаем,  автор приносит  искреннюю благодарность своему
научному руководителю Абдулла Азамовичу Азамову за постановку
задачи,  постоянное внимание к работе и ценные советы.

Саматов Бахром Таджиахматович, Наманганский государственный
университет, Е-mail:
 samatov57@inbox.ru

\end{document}